\documentclass[12pt]{amsart}

\usepackage{amsmath,amssymb,CJK}

\newtheorem{theorem}{Theorem}[section]
\newtheorem{lemma}[theorem]{Lemma}

\theoremstyle{definition}
\newtheorem{definition}[theorem]{Definition}

\theoremstyle{remark}
\newtheorem{remark}{Remark}[section]

\numberwithin{equation}{section}

\title{Quotients of Hurwitz Primes}

\author{Minghao Pan}
\address{Department of Mathematics, University of California, Los Angeles, CA 90095, United States}
\email{minghaopan@g.ucla.edu}

\author{Wentao Zhang}
\address{Shenzhen Middle School, No.18 Shenzhong Street, Luohu, Shenzhen, Guangdong, 518001, China}
\email{wtzhang@shenzhong.net}

\date{\today}

\begin{document}

\maketitle

\begin{abstract}
Quotient sets have attracted the attention of mathematicians in the past three decades. The set of quotients of primes is dense in the positive real numbers and the set of all quotients of Gaussian primes is also dense in the complex plane. Sittinger has proved that the set of quotients of primes in an imaginary quadratic ring is dense in the complex plane and the set of quotients of primes in a real quadratic number ring is dense in $\mathbb{R}.$ An interesting open question is introduced by Sittinger: Is the set of quotients of Hurwitz primes dense in the quaternions? In this paper, we answer the question and prove that the set of all quotients of Hurwitz primes is dense in the quaternions.
\end{abstract}

\section{Introduction}
Quotient sets like $\{p/q:p,q\mathrm{\ are\ primes}\}$ have attracted the attention of mathematicians in the past three decades. It has been proved (or observed) many times that the set $\{p/q:p,q\mathrm{\ are\ primes}\}$ is dense in the positive real numbers (e.g., \cite[Exercise 218]{de}, \cite[Corollary 5]{gs}, \cite[Theorem 4]{ho} ). In 2013, Garcia
\cite{ga1} considered the set of all quotients of Gaussian primes and proved that it is dense in the complex plane. Later Garcia and Luca \cite{ga2} proved that the set of quotients of nonzero Fibonacci numbers is dense in the $p$-adic numbers
for every prime $p$. Sanna \cite{sa} generalized Garcia and Luca's result and proved that for any integer $k\ge2$ and any prime number $p$, the set of quotients of nonzero
$k$-generalized Fibonacci numbers is dense in the $p$-adic numbers.

Recently, Sittinger \cite{si} proved that the set of quotients of primes in an imaginary quadratic ring is dense in the complex plane and the set of quotients of primes in a real quadratic number ring
 is dense in $\mathbb{R}.$ Sittinger also asked an interesting open question in his paper: Is the set of quotients of Hurwitz primes dense in the quaternions (see Section 2 for the definitions)? In this paper, we answer Sittinger's question and prove the following theorem.
 \begin{theorem}\label{thm1}
 The set of all quotients of Hurwitz primes is dense in the quaternions.
 \end{theorem}

\begin{remark}
As the multiplication of quaternions is not commutative, for any two non-zero quaternions $\mathfrak{a},\mathfrak{b}$, their quotient could be defined by  $\frac{\mathfrak{b}}{\mathfrak{a}}=\mathfrak{b}\mathfrak{a}^{-1}$ or $\frac{\mathfrak{b}}{\mathfrak{a}}=\mathfrak{a}^{-1}\mathfrak{b}$ and Theorem \ref{thm1} holds for both cases. In this paper, we only prove the former case and the proof of the latter case is very similar with obvious modifications.
\end{remark}

\begin{remark}
In fact we prove a slightly stronger result than Theorem \ref{thm1}. Hurwitz quaternions could be divided into two disjoint subsets (see Section 2 for the notation and definitions)
$$
H_1=\left\{{x_1+x_2i+x_3j+x_4k}:x_1,x_2,x_3,x_4\in\mathbb{Z}\right\}
$$
and
$$
H_2=\left\{{x_1+x_2i+x_3j+x_4k}:x_1,x_2,x_3,x_4\in\mathbb{Z}+\frac12\right\}.
$$
Our proof indicates that the set
$$
\left\{\frac{\mathfrak{p}}{\mathfrak{q}}:{\mathfrak{p}}\ \textrm{and}\ {\mathfrak{q}}\ \textrm{are Hurwitz primes in } H_1\right\}
$$
is dense in the quaternions. Moreover, for any two Hurwitz primes ${\mathfrak{p}},{\mathfrak{q}}\in H_1$ with odd norms, it is easy to see that $\mathfrak{pu},\mathfrak{qu}$ are Hurwitz primes belonging to $H_2$ and $\frac{\mathfrak{pu}}{\mathfrak{qu}}=(\mathfrak{pu})(\mathfrak{qu})^{-1}=\frac{\mathfrak{p}}{\mathfrak{q}}$, where  $\mathfrak{u}=\frac{1+i+j+k}{2}$ is a unit in Hurwitz quaternions. Therefore, the set
$$
\left\{\frac{\mathfrak{p}}{\mathfrak{q}}:{\mathfrak{p}}\ \textrm{and}\ {\mathfrak{q}}\ \textrm{are Hurwitz primes in } H_2\right\}
$$
is also dense in the quaternions.
\end{remark}

\section{Hurwitz quaternions}
In this section, we introduce some properties of quaternions and most of the materials can be found in \cite{co}.

The quaternions were discovered by Irish mathematician Hamilton in 1843. They have been widely used in the electrodynamics, general relativity, navigation, satellite attitude control and other fields.

\begin{definition}
The set of quaternions is defined as
$$
Q=\left\{x_1+x_2i+x_3j+x_4k:x_1,x_2,x_3,x_4\in\mathbb{R}\right\}
$$
where $i,j,k$ commute with every real number and satisfy
\begin{equation*}\label{ijk}
ijk=i^2=j^2=k^2=-1.
\end{equation*}
\end{definition}
Let $\mathfrak{a}=a_1+a_2i+a_3j+a_4k$ and $\mathfrak{b}=b_1+b_2i+b_3j+b_4k$ be any two quaternions. The addition of quaternions is defined by
$$
\mathfrak{a}+\mathfrak{b}=a_1+b_1+(a_2+b_2)i+(a_3+b_3)j+(a_4+b_4)k.
$$
For any real number $\lambda$, the scalar multiplication is defined by
$$
\lambda\mathfrak{a}=\lambda a_1+\lambda a_2i+\lambda a_3j+\lambda a_4k.
$$
Then the quaternions form a vector space with these two operations. Moreover, we can define the multiplication of quaternions by
\begin{align*}
\mathfrak{a}\mathfrak{b}
=&(a_1b_1-a_2b_2-a_3b_3-a_4b_4)+(a_1b_2+a_2b_1+a_3b_4-a_4b_3)i\\
&+(a_1b_3+a_3b_1+a_4b_2-a_2b_4)j+(a_1b_4+a_4b_1+a_2b_3-a_3b_2)k.
\end{align*}
Clearly, we have
$$
ij=k\mathrm{\ and\ }ji=-k
$$
so the multiplication of quaternions is not commutative.

For any $\mathfrak{a}=a_1+a_2i+a_3j+a_4k\in Q$, $\overline{\mathfrak{a}}=a_1-a_2i-a_3j-a_4k$ is called the conjugate of $\mathfrak{a}$. It is easy to see that
$$
\mathfrak{a}\overline{\mathfrak{a}}=a_1^2+a_2^2+a_3^2+a_4^2.
$$
\begin{definition}
For any $\mathfrak{a}=a_1+a_2i+a_3j+a_4k\in Q$, its norm is defined by
$$
\|\mathfrak{a}\|=\mathfrak{a}\overline{\mathfrak{a}}=a_1^2+a_2^2+a_3^2+a_4^2.
$$
\end{definition}
The norm induces a metric $d(\mathfrak{a},\mathfrak{b})=|\mathfrak{a}-\mathfrak{b}|$ on the quaternions by
$$
|\mathfrak{a}-\mathfrak{b}|=\sqrt{\|\mathfrak{a}-\mathfrak{b}\|}
$$
and the quaternions form a metric space.
\begin{definition}
A subset $D$ of quaternions is said to be dense in the quaternions if for any quaternion $\mathfrak{a}$ and any $\varepsilon>0$, there exists a quaternion $\mathfrak{b}\in D$ such that
$$
|\mathfrak{a}-\mathfrak{b}|<\varepsilon.
$$
\end{definition}
\begin{definition}
For any $\mathfrak{a}=a_1+a_2i+a_3j+a_4k\in Q$ and ${\|\mathfrak{a}\|}\ne 0$, its inverse is defined by
$$
\mathfrak{a}^{-1}=\frac{\overline{\mathfrak{a}}}{\|\mathfrak{a}\|}=\frac{a_1-a_2i-a_3j-a_4k}{\|\mathfrak{a}\|}.
$$
\end{definition}
In this paper, the quotient of two quaternions is defined by $$\frac{\mathfrak{b}}{\mathfrak{a}}=\mathfrak{b}\mathfrak{a}^{-1}.$$

One interesting subset of quaternions is the set of Hurwitz quaternions which was introduced by Hurwitz in 1919.
\begin{definition}
The set of Hurwitz quaternions $H$ is a subset of quaternions, defined as
$$
H=\left\{{x_1+x_2i+x_3j+x_4k}:x_1,x_2,x_3,x_4\in\mathbb{Z}\mathrm{\ or\ }x_1,x_2,x_3,x_4\in\mathbb{Z}+\frac12\right\}.
$$
We say that $\mathfrak{a}$ is a unit in $H$ if $\|\mathfrak{a}\|=1$.
\end{definition}
It is easy to see that for any $\mathfrak{a}\in H$, $\|\mathfrak{a}\|\in\mathbb{Z}$. Moreover, we have the following result.
\begin{lemma}\label{lemma 1}
Let $n$ be any positive integer. Then the number of Hurwitz quaternions with norm $n$ is $$24\sum\limits_{d|n\atop 2\nmid d}d.$$
\end{lemma}
\begin{definition}
We say that $\mathfrak{p}\in H$ is a Hurwitz prime if $\mathfrak{p}$ is not zero or a unit and is not a product of non-units in $H$.
\end{definition}
We have the following result to determine whether a Hurwitz quaternion is a Hurwitz prime.
\begin{lemma}\label{lemma 2}
For any $\mathfrak{p}\in H$, $\mathfrak{p}$ is a Hurwitz prime if and only if $\|\mathfrak{p}\|$ is a prime number.
\end{lemma}

\section{Preliminaries}
In this section, we introduce some tools which will be used later. We begin with some well-known properties of $\mathbb{R}^4$.
For any two vectors $\overrightarrow{x}=(x_1,x_2,x_3,x_4),\overrightarrow{y}=(y_1,y_2,y_3,y_4)\in\mathbb{R}^4$, the metric is defined by
$$
|\overrightarrow{x}-\overrightarrow{y}|
=\sqrt{(x_1-y_1)^2+(x_2-y_2)^2+(x_3-y_3)^2+(x_4-y_4)^2}
$$
and the inner product is defined by
$$
\langle\overrightarrow{x},\overrightarrow{y}\rangle=x_1y_1+x_2y_2+x_3y_3+x_4y_4.
$$
Clearly $|\overrightarrow{x}|=\sqrt{\langle\overrightarrow{x},\overrightarrow{x}\rangle}$.
It is well-known that
\begin{equation}\label{para}
|\overrightarrow{x}-\overrightarrow{y}|^2
=|\overrightarrow{x}|^2+|\overrightarrow{y}|^2
-2
{\langle\overrightarrow{x},\overrightarrow{y}\rangle}
.
\end{equation}

Define a map $\sigma$ from $Q$ to $\mathbb{R}^4$ by
\begin{align*}
\sigma:\qquad\qquad Q\qquad&\rightarrow\qquad\mathbb{R}^4\\
x_1+x_2i+x_3j+x_4k\ &\rightarrow\ (x_1,x_2,x_3,x_4).
\end{align*}
Then it is easy to see that $\sigma$ is an isomorphism. Moreover, $\sigma$ is also an isometry and
\begin{equation}\label{equiv}
|\sigma(\mathfrak{a})-\sigma(\mathfrak{b})|=|\mathfrak{a}-\mathfrak{b}|
\end{equation}
for any $\mathfrak{a},\mathfrak{b}\in Q$.

Next, we introduce our main tool.
Denote by $S$ the four dimensional hypersphere
$$
x_1^2+x_2^2+x_3^2+x_4^2=1.
$$
For any $0<\theta<\pi$ and $\overrightarrow{x}\in\mathbb{R}^4$, define
$$
\Omega(\overrightarrow{x},\theta)
=\left\{\overrightarrow{y}\in\mathbb{R}^4:|\overrightarrow{y}|=1\mathrm{\ and\ }\arccos
\frac{\langle \overrightarrow{x},\overrightarrow{y}\rangle}
{|\overrightarrow{x}||\overrightarrow{y}|}\le \theta\right\}.
$$
$\Omega(\overrightarrow{x},\theta)$ is a hyperspherical cap in $S$ and denote by $A(\Omega(\overrightarrow{x},\theta))$ its surface area. Clearly $A(\Omega(\overrightarrow{x},\theta))$ is a positive real number and only depends on $\theta$ and $\overrightarrow{x}$.

Define
\begin{equation}\label{def r omega}
r(n,\Omega(\overrightarrow{x},\theta))
=\#\left\{\overrightarrow{y}\in\mathbb{Z}^4:|\overrightarrow{y}|=\sqrt{n}
\mathrm{\ and\ }\frac{\overrightarrow{y}}{\sqrt{n}}\in\Omega(\overrightarrow{x},\theta)\right\}.
\end{equation}
The following theorem is a special case of \cite[Theorem 1]{fo} with $Q(X)=x_1^2+x_2^2+x_3^2+x_4^2$ and $\Omega=\Omega(\overrightarrow{x},\theta)$.
\begin{theorem}\label{thm fo}
Let notation be as above. For any positive integer $n$ with $(n,2)=1$ and $\varepsilon>0$, we have
$$
r(n,\Omega(\overrightarrow{x},\theta))
=r(n)\frac{A(\Omega(\overrightarrow{x},\theta))}{A(S)}\left(1+O\left(n^{-1/7+\varepsilon}\right)\right),
$$
where $r(n)$ is the number of integral solutions of $x_1^2+x_2^2+x_3^2+x_4^2=n$ and $A(S)$ is the surface area of $S$.
\end{theorem}
\begin{remark}
By the famous Jacobi's four-square theorem, we have
\begin{equation}\label{jacobi}
r(n)=\left\{ \begin{aligned}
         &8\sum_{m|n}m \ \ \ \mathrm{if\ }n\ \mathrm{is\ odd},&\\
           & 24\sum_{m|n\atop2\nmid m}m\ \ \ \mathrm{if\ }n\ \mathrm{is\ even}.&
                          \end{aligned} \right.
\end{equation}
\end{remark}

$ \\ $

\section{Proof of Theorem \ref{thm1}}
It is sufficient to prove that for any quaternion $\mathfrak{h}$ and any $\varepsilon>0$, there exist two Hurwitz primes $\mathfrak{p},\mathfrak{q}$ such that $$\left|\mathfrak{h}-\frac{\mathfrak{p}}{\mathfrak{q}}\right|<\varepsilon.$$

We first consider the case $\|\mathfrak{h}\|=0$. Since the set of all quotients of prime numbers is dense in positive real numbers, there exist two prime numbers $p,q$ such that $p/q<\varepsilon^2$. By Lemma \ref{lemma 1} and Lemma \ref{lemma 2}, there exist two Hurwitz primes $\mathfrak{p},\mathfrak{q}$ such that $\|\mathfrak{p}\|=p$, $\|\mathfrak{q}\|=q$ and $$\left|\frac{\mathfrak{p}}{\mathfrak{q}}\right|=\sqrt{\frac pq}<\varepsilon.$$

In what follows, we assume $\|\mathfrak{h}\|\neq0$ and
without loss of generality, we assume
\begin{equation}\label{varepsilon}
\varepsilon<\min(\|\mathfrak{h}\|,1/\|\mathfrak{h}\|)\le1.
\end{equation}
Put
\begin{equation}\label{varepsilon1}
\varepsilon_1=\frac{\varepsilon^2}{10(\|\mathfrak{h}\|+\varepsilon)}
\le 1.
\end{equation}
By Theorem \ref{thm fo} and \eqref{jacobi}, for any positive odd integer $n$, we have
$$
r(n,\Omega(\sigma(\mathfrak{h}),\varepsilon_1))
=8\frac{A(\Omega(\sigma(\mathfrak{h}),\varepsilon_1))}{A(S)}\sum_{m|n}m
\left(1+O\left(n^{-1/7+\varepsilon_1}\right)\right).
$$
Since $\frac{A(\Omega(\sigma(\mathfrak{h}),\varepsilon_1))}{A(S)}$ is positive and only depends on $\varepsilon$ and $\mathfrak{h}$, there exists $N_1=N_1(\varepsilon,\mathfrak{h})$ such that
$$
r(n,\Omega(\sigma(\mathfrak{h}),\varepsilon_1))>1
$$
if $n>N_1$.
By similar arguments, there exists $N_2=N_2(\varepsilon_1,\overrightarrow{e_1})$ such that
$$
r(n,\Omega(\overrightarrow{e_1},\varepsilon_1))>1
$$
if $n>N_2$, where $\overrightarrow{e_1}=(1,0,0,0)$.
Moreover, by the Prime Number Theorem, there exists $N_3=N_3(\varepsilon,\mathfrak{h})$ such that the interval $$(n(\|\mathfrak{h}\|-\varepsilon^2/10),n(\|\mathfrak{h}\|+\varepsilon^2/10))$$ contains at least one prime number if $n>N_3$.

Let $q$ be a prime number satisfying $$q>\max\left(\frac{N_1}{\|\mathfrak{h}\|-\varepsilon^2/10},N_2,N_3\right).$$
Then
\begin{equation}\label{existence q}
r(q,\Omega(\overrightarrow{e_1},\varepsilon_1))>1
\end{equation}
and there exists a prime
\begin{equation}\label{p}
p\in (q(\|\mathfrak{h}\|-\varepsilon^2/10),q(\|\mathfrak{h}\|+\varepsilon^2/10)).
\end{equation}
By our choice of $q$, we get that
$$
p>q(\|\mathfrak{h}\|-\varepsilon^2/10)>N_1.
$$
Hence, we obtain
\begin{equation}\label{existence p}
r(p,\Omega(\sigma(\mathfrak{h}),\varepsilon_1))>1.
\end{equation}

By \eqref{def r omega}, \eqref{existence q} and \eqref{existence p}, there exist $$\overrightarrow{x}=(x_1,x_2,x_3,x_4)\in\mathbb{Z}^4\mathrm{\ and\ }\overrightarrow{y}=(y_1,y_2,y_3,y_4)\in\mathbb{Z}^4$$ such that $|\overrightarrow{x}|=\sqrt{q}$, $|\overrightarrow{y}|=\sqrt{p}$,
\begin{equation}\label{upper bound 1}
\arccos\frac{x_1}{\sqrt{q}}
=
\arccos\frac{\langle \overrightarrow{x},\overrightarrow{e_1}\rangle}{|\overrightarrow{x}|}
\le \varepsilon_1
\end{equation}
and
\begin{equation}\label{upper bound 2}
\arccos
\frac{\langle \overrightarrow{y},\sigma(\mathfrak{h})\rangle}
{|\overrightarrow{y}||\sigma(\mathfrak{h})|}\le \varepsilon_1.
\end{equation}
By \eqref{upper bound 1}, we have
\begin{align}\label{upper bound 3}
0 \le 1-\frac{x_1}{\sqrt{q}}\le 1-\cos\varepsilon_1=2\sin^2\frac{\varepsilon_1}{2}
\le\frac{\varepsilon_1^2}{2}
\end{align}
and for $\ell=2,3,4$
\begin{align}\label{upper bound 4}
0 \le \frac{x_\ell^2}{{q}}\le 1-\frac{x_1^2}{{q}}\le 1-\cos^2\varepsilon_1
=\sin^2\varepsilon_1\le\varepsilon_1^2.
\end{align}
Here we have used the well-known inequality $0\le \sin t\le t$ if $0\le t\le 1$.
Moreover, by \eqref{p}, we have
\begin{equation}\label{norm 1}
\left|\frac{\overrightarrow{y}}{\sqrt{q}}\right|
=\sqrt{\frac pq}\le \sqrt{\|\mathfrak{h}\|+\varepsilon^2/10}
\end{equation}
and
\begin{align*}
\left(|\sigma(\mathfrak{h})|-\left|\frac{\overrightarrow{y}}{\sqrt{q}}\right|\right)^2
=\left(\frac{|\sigma(\mathfrak{h})|^2-\left|\frac{\overrightarrow{y}}{\sqrt{q}}\right|^2}
{|\sigma(\mathfrak{h})|+\left|\frac{\overrightarrow{y}}{\sqrt{q}}\right|}\right)^2
\le\left(\frac{\|\mathfrak{h}\|-\frac{p}{{q}}}
{|\sigma(\mathfrak{h})|}\right)^2\le \frac{\varepsilon^4}{100\|\mathfrak{h}\|}.
\end{align*}
Therefore, by \eqref{para}, \eqref{upper bound 2} and the last inequality in \eqref{upper bound 3}, we obtain
\begin{align}\label{difference 1}
\left|\sigma(\mathfrak{h})-\frac{\overrightarrow{y}}{\sqrt{q}}\right|^2
&=|\sigma(\mathfrak{h})|^2
+\left|\frac{\overrightarrow{y}}{\sqrt{q}}\right|^2
-2|\sigma(\mathfrak{h})|\left|\frac{\overrightarrow{y}}{\sqrt{q}}\right|
\frac{\langle\sigma(\mathfrak{h}),\frac{\overrightarrow{y}}{\sqrt{q}}\rangle}
{|\sigma(\mathfrak{h})|\left|\frac{\overrightarrow{y}}{\sqrt{q}}\right|}\nonumber\\
&\le|\sigma(\mathfrak{h})|^2
+\left|\frac{\overrightarrow{y}}{\sqrt{q}}\right|^2
-2|\sigma(\mathfrak{h})|\left|\frac{\overrightarrow{y}}{\sqrt{q}}\right|\cos\varepsilon_1\nonumber\\
&=\left(|\sigma(\mathfrak{h})|-
\left|\frac{\overrightarrow{y}}{\sqrt{q}}\right|\right)^2
+2|\sigma(\mathfrak{h})|\left|\frac{\overrightarrow{y}}{\sqrt{q}}\right|(1-\cos\varepsilon_1)\nonumber\\
&\le\frac{\varepsilon^4}{100\|\mathfrak{h}\|}
+{\varepsilon_1^2}\sqrt{\|\mathfrak{h}\|(\|\mathfrak{h}\|+\varepsilon^2/10)}
\le \frac{\varepsilon^4}{50\|\mathfrak{h}\|}\le \frac{\varepsilon^2}{9}.
\end{align}
Here we have applied \eqref{varepsilon} and \eqref{varepsilon1} in the last two steps.

Put
$$
\mathfrak{q}=x_1+x_2i+x_3j+x_4k
$$
and
$$
\mathfrak{p}=y_1+y_2i+y_3j+y_4k.
$$
Then $\|\mathfrak{p}\|=p$ and $\|\mathfrak{q}\|=q$. By Lemma \ref{lemma 2}, $\mathfrak{p}$ and $\mathfrak{q}$ are Hurwitz primes. Furthermore, by the triangle inequality we have
\begin{align}\label{upper bound 5}
\left|\mathfrak{h}-\frac{\mathfrak{p}}{\mathfrak{q}}\right|
&=\left|\mathfrak{h}-\frac{\mathfrak{p}(x_1-x_2i-x_3j-x_4k)}{\|\mathfrak{q}\|}\right|\nonumber\\
&\le \left|\mathfrak{h}-\frac{x_1}{q}\mathfrak{p}\right|
+\left|\frac{\mathfrak{p}(x_2i)}{q}\right|
+\left|\frac{\mathfrak{p}(x_3j)}{q}\right|
+\left|\frac{\mathfrak{p}(x_4k)}{q}\right|.
\end{align}
By \eqref{upper bound 4} and \eqref{p}, we obtain
\begin{align}\label{upper bound 6}
&\left|\frac{\mathfrak{p}(x_2i)}{q}\right|
+\left|\frac{\mathfrak{p}(x_3j)}{q}\right|
+\left|\frac{\mathfrak{p}(x_4k)}{q}\right|\nonumber\\
&=\sum_{\ell=2}^4\sqrt{\frac{x_\ell^2\|\mathfrak{p}\|}{q^2}}
=\sqrt{\frac{p}{q}}\sum_{\ell=2}^4\sqrt{\frac{x_\ell^2}{q}}\nonumber\\
&\le 3\sqrt{(\|\mathfrak{h}\|+\varepsilon^2/10)}\varepsilon_1\le \frac{3\varepsilon^2}{10\sqrt{\|\mathfrak{h}\|}}\le\frac{\varepsilon}{3}.
\end{align}
Here we have applied \eqref{varepsilon} and \eqref{varepsilon1} in the last two steps.
On the other hand, by \eqref{equiv}, \eqref{difference 1}, \eqref{upper bound 3} and \eqref{norm 1}, we get
\begin{align}\label{upper bound 7}
\left|\mathfrak{h}-\frac{x_1}{q}\mathfrak{p}\right|
=\left|\sigma(\mathfrak{h})-\frac{x_1}{q}\sigma(\mathfrak{p})\right|
&=\left|\sigma(\mathfrak{h})-\frac{\overrightarrow{y}}{\sqrt{q}}
+\left(1-\frac{x_1}{\sqrt{q}}\right)\frac{\overrightarrow{y}}{\sqrt{q}}\right|
\nonumber\\
&\le\left|\sigma(\mathfrak{h})-\frac{\overrightarrow{y}}{\sqrt{q}}
\right|
+\left(1-\frac{x_1}{\sqrt{q}}\right)\left|\frac{\overrightarrow{y}}{\sqrt{q}}\right|\nonumber\\
&\le \frac{\varepsilon}{3}+\frac{\varepsilon_1^2}{2}\sqrt{(\|\mathfrak{h}\|+\varepsilon^2/10)}
\le\frac{2\varepsilon}{3}.
\end{align}
Here we have applied \eqref{varepsilon} and \eqref{varepsilon1} again in the last one step.
Combining \eqref{upper bound 5}, \eqref{upper bound 6} and \eqref{upper bound 7}, we have
$$
\left|\mathfrak{h}-\frac{\mathfrak{p}}{\mathfrak{q}}\right|\le \frac{2\varepsilon}{3}+\frac{\varepsilon}{3}=\varepsilon.
$$
The proof is complete.

\section{Acknowledgement}
It is our pleasure to thank professor Yingnan Wang who is from Shenzhen Univeristy  for his help and helpful advice throughout this project.


\bibliography{Biblio2}
\bibliographystyle{alpha}

\end{document}